# General Fractional Calculus: Multi-Kernel Approach

**Vasily E. Tarasov** [1, 2,*]

[1] Skobeltsyn Institute of Nuclear Physics, Lomonosov Moscow State University, Moscow 119991, Russia; tarasov@theory.sinp.msu.ru; Tel.: +7-495-939-5989

[2] Faculty "Information Technologies and Applied Mathematics", Moscow Aviation Institute, National Research University, Moscow 125993, Russia

**Abstract:** For the first time, a general fractional calculus of arbitrary order was proposed by Yuri Luchko in the works *Mathematics 9(6) (2021) 594* and *Symmetry 13(5) (2021) 755*. In these works, the proposed approaches to formulate this calculus are based either on the power of one Sonine kernel, or the convolution of one Sonine kernel with the kernels of the integer-order integrals. To apply general fractional calculus, it is useful to have a wider range of operators, for example, by using the Laplace convolution of different types of kernels. In this paper, an extended formulation of the general fractional calculus of arbitrary order is proposed. Extension is achieved by using different types (subsets) of pairs of operator kernels in definitions general fractional integrals and derivatives. For this, the definition of the Luchko pair of kernels is somewhat broadened, which leads to the symmetry of the definition of the Luchko pair. The proposed set of kernel pairs are subsets of the Luchko set of kernel pairs. The fundamental theorems for the proposed general fractional derivatives and integrals are proved.

## 1. Introduction

The theory of integro-differential operators and equations is important tool to describe systems and processes with non-locality in space and time. Among such operators, an important role is played by the integrals and derivatives of non-integer orders [1, 2, 3, 4, 5], [6, 7]. These operators are called fractional derivatives (FD) and fractional integrals (FI). For these operators, generalizations of first and second fundamental theorems of standard calculus are satisfied. This is one of the main reasons to state that these operators form a calculus, which is called fractional calculus.

Equations with derivatives of non-integer orders with respect to time and space are important tools to describe non-locality in time and space in physics [8, 9], biology [10], and economics [11, 12], for example. To describe various physical, biological, economic phenomena with nonlocality, it is important to have a wide range of operators that allow us to describe various types of nonlocality [13]. Nonlocality is determined by the form of the kernel of the operator, which are fractional integrals (FI) and fractional derivatives (FD) of non-integer orders. Therefore it is important to have a fractional calculus that allows us to describe non-locality in a general form.

An important turning point in formulation of general fractional calculus (GFC) was the results obtained by Anatoly N. Kochubei in his work [14] (see also [15]) in 2011. In this work, the general fractional integral (GFI) and general fractional derivatives (GFD) of the Riemann-Liouville and Caputo type are defined, for which the general fundamental theorems are proved (see Theorem 1 in [14]). In addition, the relaxation and diffusion equations [14, 16], and then the growth equation [17], which are contain GFD, are solved. In fact, the term "general fractional calculus" (GFC) was introduced in article [14]. This approach to GFC is based on the concept of Sonine pairs of mutually associated kernels proposed in the work [18] (see also [19]). The integral equations of the first kind with Sonine kernels, and the GFI and GFD of the Liouville and Marchaud type are described in [20, 21]. After Kochubei article [14], works on general fractional calculus and some of its applications began to be published (see [22, 23, 24, 25, 26] and references therein).

The next revolutionary step in the construction of the general fractional calculus was proposed in the works of Yuri Luchko in 2021 [27, 28, 29]. In artiles [27, 28] a general fractional integrals and derivatives of arbitrary order have been proposed. The general fundamental theorems of GFD are proved for the GFI and GFDs of Riemann-Liouville and Caputo in [27, 28]. Operational calculus for equations with general fractional derivatives with the Sonine kernels is proposed in [29].

In articles [27, 28], two possible approaches to construct general fractional integrals and derivatives of arbitrary order, which satisfy general fundamental theorems of GFC, are proposed. These approaches are based on building the Luchko pairs $(M(t), N(t)) \in L_n$ with $n > 1$ for the GFI and GFD kernels from the Sonine pairs of kernels $(\mu(t), \nu(t)) \in S_{-1} = L_1$, [27, 28]. These two possible approaches to define kernels of GDI and GFD, which are proposed in [27, 28], can be briefly described as follows.

1) In article [27], the kernels $M(t), N(t)$ of GFI and GFD, which belong to the Luchko set $L_n$, are considered in the form

$$M(t) = (\mu_1 * \ldots * \mu_n)(t), \quad N(t) = (\nu_1 * \ldots * \nu_n)(t), \tag{1}$$

for the case $\mu_k(t) = \mu(t)$ and $\nu_k(t) = \nu(t)$ for all $k = 1, \ldots, n$ such that

$$\mu(t), \nu(t) \in C_{-1,0}(0, \infty), \quad \text{and} \quad (\mu * \nu)(t) = \{1\}, \tag{2}$$

where $*$ is the Laplace convolution, $\{1\}$ is the function that is equal to 1 for $t \geq 0$, and the function $f(t)$ belongs to the space $C_{-1,0}(0, \infty)$, if it can be represented as $f(t) = t^p g(t)$, where $-1 < p < 0$ and $g(t) \in C[0, \infty)$ for $t > 0$.

2) In article [28], the kernels $M(t), N(t)$ of GFI and GFD, which belong to the Luchko set $L_n$, are considered in the form

$$M(t) = (\{1\}^{n-1} * \mu)(t), \quad \text{and} \quad N(t) = \nu(t), \tag{3}$$

where

$$\mu(t), \nu(t) \in C_{-1,0}(0, \infty), \quad \text{and} \quad (\mu * \nu)(t) = \{1\}^n = \frac{t^{n-1}}{(n-1)!} \tag{4}$$

3) In [28], the set $L_n$, of kernel pairs $(M(t), N(t))$ is defined by the conditions

$$M(t) \in C_{-1}(0, \infty), \quad N(t) \in C_{-1,0}(0, \infty), \quad \text{and} \quad (\mu * \nu)(t) = \{1\}^n. \tag{5}$$

The requirement $N(t) \in C_{-1,0}(0, \infty)$ instead of $N(t) \in C_{-1}(0, \infty)$ leads to the fact that $(\mu^n(t), \nu^n(t)) \in L_n$ only if $N(t) = \nu^n(t) \in C_{-1}(0, \infty)$, which is very strong and restrictive condition. The function $f(t)$ belongs to the space $C_{-1}(0, \infty)$, if it can be represented in the form $f(t) = t^p g(t)$, where $p > -1$ and $g(t) \in C[0, \infty)$ for $t > 0$.

In works [27, 28], the proposed approaches are based either on the powers of one Sonine kernel [27], or the convolution of one Sonine kernel with the kernels of the integer-order integrals [28]. In applications of fractional calculus and GFC, it is useful to have a wider range of operators, for example, by using the Laplace convolution of different types of kernels.

Let us give some examples of possible expansions of operator kernels of the general fractional calculus. In the beginning we proposed to define the Luchko pair of kernels is somewhat broadened, which leads to the symmetry of the definition of the Luchko pairs, by using $M(t), N(t) \in C_{-1}(0, \infty)$ in Definition 5. Then we can consider, for example, the following pairs of kernels from the Luchko set $L_n$.

I) As an example of a generalization of the first approach, we can use the kernels that are the Laplace convolution of different types of kernels. For example, one of the ways to define the kernels $(M(t), N(t)) \in \mathcal{L}_n$ with $n > 1$ is to remove the restrictions $\mu_k(t) = \mu(t)$ and $\nu_k(t) = \nu(t)$ for all $k = 1, \ldots, n$, which are used in (1).

II) Another example is removing the restriction on using only one pair of Sonine kernel in (3). For example, we can consider the Luchko pairs $(M(t), N(t)) \in L_n$ with $n > 1$ in the form

$$M(t) = (\{1\}^{n-m} * \mu_1 * \ldots * \mu_m)(t) \quad \text{and} \quad N(t) = (\nu_1 * \ldots * \nu_m)(t) \tag{6}$$

instead of (3), where
$$\mu_k(t), \nu_k(t) \in C_{-1,0}(0,\infty), \quad \text{and} \quad (\mu_k * \nu_k)(t) = \{1\} \quad \text{for all} \quad k = 1, \dots, m. \tag{7}$$

III) As a more general example of the pair of kernels from $L_n$, we can consider the Laplace convolutions
$$M(t) = (M_{k_1} * \dots * M_{k_p})(t), \quad N(t) = (N_{k_1} * \dots * N_{k_p})(t), \tag{8}$$
of kernel pairs $(M_{k_j}(t), N_{k_j}(t))$ from the Luchko sets $L_{k_j}$ such that
$$M_{k_j}(t), N_{k_j}(t) \in C_{-1}(0,\infty) \quad \text{and} \quad \left(M_{k_j}(t) * N_{k_j}\right)(t) = \{1\}^{k_j} \quad \text{where} \quad \sum_{j=1}^{p} k_j = n \tag{9}$$
for all $j = 1, \dots p$.

These possible approaches of extensions are based on the statement: the triple $\mathcal{R}_{-1} = (C_{-1}(0,\infty), *, +)$, is a commutative ring without divisors of zero [30, 27], where the multiplication $*$ is the Laplace convolution and $+$ the standard addition of functions. These examples and other possible approaches to expanding the variety of types of kernels of operators of general fractional calculus and, thus, nonlocality, are important for describing systems and processes with nonlocality in space and time.

In this paper, an extended formulation of the general fractional calculus of arbitrary order is proposed as an extension of the Luchko approaches, which is described in [27, 28]. Extension is achieved by using different types (subsets) of Sonine and Luchko pairs of kernels in definitions general fractional integrals and derivatives of multi-kernel form. For this, the definition of the Luchko pair of kernels is somewhat broadened, which leads to the symmetry of the definition of the Luchko pair. The proposed sets of kernel pairs are subsets of the Luchko set of kernel pairs. The first and second fundamental theorems for the proposed general fractional derivatives and integrals of multi-kernel form are proved in this paper.

## 2 Luchko Set of Kernel Pairs and its Subsets

Let us give definitions of some extеansion of the concept of the Luchko pairs of kernels $L_n$, which is given in [28] and some subset of kernels from $L_n$. In these definitions, we use the function spaces $C_{-1,0}(0,\infty)$ and $C_{-1}(0,\infty)$. The function $f(t)$ belongs to the space $C_{-1,0}(0,\infty)$ or $C_{-1}(0,\infty)$, if it can be represented in the form $f(t) = t^p g(t)$, where $g(t) \in C[0,\infty)$ for $t > 0$ and $-1 < p < 0$ or $p > -1$, respectively.

**Definition 1** *Let $M(t), N(t) \in C_{-1}(0,\infty)$ and the Luchko condition*
$$(M * N)(t) = \{1\}^n \tag{10}$$
*holds for $t \in (0,\infty)$, where*
$$\{1\}^n = h_n(t) = \frac{t^{n-1}}{(n-1)!}, \quad (M * N)(t) = \int_0^t M(t-\tau) N(\tau) \, d\tau. \tag{11}$$
*Then the set of such kernel pairs $(M(t), N(t))$ will be called the Luchko pairs and it will be denoted as $L_n$. The kernels $M(t), N(t)$ will be called the Luchko kernels, and the pair $(M(t), N(t))$ will be called the Luchko pairs.*

**Remark 1** *Note that Definition 1 is a minor generalization of the definition of the set of kernel pairs $L_n$, which is proposed in [28, p.7]. In article [28, p.7], the set $L_n$ is defined as a set of the kernels $M(t), \in C_{-1}(0,\infty)$ and $N(t) \in C_{-1,0}(0,\infty) \subset C_{-1}(0,\infty)$ that satisfy condition (10).*

**Remark 2** *Note that the pairs of kernels $\mu(t), \nu(t) \in C_{-1,0}(0,\infty)$ that satisfy condition (10) with $n = 1$ is called the Sonine pairs and the set of such pairs is denoted as $S_{-1}$. Using the property $C_{-1,0}(0,\infty) \subset C_{-1}(0,\infty)$, the pairs $(\mu(t), \nu(t))$ belong to the Luchko set $L_1$.*

Let us consider the following approach to construct pairs $(M(t), N(t)) \in L_n$ with $n > 1$ by using the Sonine pair of the kernels $(\mu_k(t), \nu_k(t)) \in L_1$, where $k \in \mathbb{N}$.

**Theorem 1** *Let $(\mu_k(t), \nu_k(t))$ be Sonine pairs of kernels from $L_1$ for all $k = 1, \dots, n$. Then the pair of kernels*
$$M_n(t) = (\mu_1 * \dots * \mu_n)(t), \quad N_n(t) = (\nu_1 * \dots * \nu_n)(t) \tag{12}$$
*belong to the Luchko set $L_n$.*

**Proof**. Using the equalities
$$(M_n * N_n)(t) = ((\mu_1 * \ldots * \mu_n) * (\nu_1 * \ldots * \nu_n))(t) =$$
$$(\mu_1 * \ldots * \mu_n * \nu_1 * \ldots * \nu_n)(t) =$$
$$((\mu_1 * \nu_1) * \ldots * (\mu_n * \nu_n))(t) =$$
$$(\{1\} * \ldots * \{1\})(t) = \{1\}^n, \tag{13}$$

we get that kernels (12) satisfy the Luchko condition (1).

Then we use that the triple $\mathcal{R}_{-1} = (C_{-1}(0, \infty), *, +)$, where the multiplication $*$ is the Laplace convolution and $+$ the standard addition of functions, is a commutative ring without divisors of zero [30, 27]. Therefore, using that $\mu_k(t), \nu_k(t) \in C_{-1}(0, \infty)$ for all $k = 1, \ldots, n$, we obtain that kernels (12) belong to $C_{-1}(0, \infty)$ and the pair $(M_n(t), N_n(t))$ belongs to $L_n$.

□

**Definition 2** *The set of kernels pairs $(M(t), N(t))$, which belong to $L_n$ and can be represented in the form*
$$M(t) = M_n(t) = (\mu_1 * \ldots * \mu_n)(t), \quad N(t) = N_n(t) = (\nu_1 * \ldots * \nu_n)(t), \tag{14}$$
*where $(\mu_k(t), \nu_k(t))$ are the Sonine pairs of kernels from $L_1$ for all $k = 1, \ldots, n$, will be denoted as $T_n$. The set $T_n$ is a subset of the Luchko set $L_n$.*

**Example 1** Let us consider the kernels
$$M_2(t) = (\mu_1 * \mu_1)(t), \quad N_2(t) = (\nu_1 * \nu_2)(t), \tag{15}$$
which belong to $T_2$, where the kernel pairs are
$$\mu_1(t) = h_\alpha(t) = \frac{t^{\alpha-1}}{\Gamma(\alpha)}, \quad \nu_1(t) = h_{1-\alpha}(t), \tag{16}$$
$$\mu_2(t) = h_\beta(t) = \frac{t^{\beta-1}}{\Gamma(\beta)}, \quad \nu_2(t) = h_{1-\beta}(t), \tag{17}$$
with $\alpha, \beta \in (0,1)$. In this case, the kernels have the form
$$M_2(t) = h_{\alpha+\beta}(t) = \frac{t^{\alpha+\beta-1}}{\Gamma(\alpha+\beta)}, \quad N_2(t) = h_{2-\alpha-\beta} = \frac{t^{1-\alpha-\beta}}{\Gamma(2-\alpha-\beta)}. \tag{18}$$

**Remark 3** *Using example 1, we can see that the pairs $(M_2(t), N_2(t))$ always belong to the Luchko set $L_2$ in the sense of Definition 1. Note that this is not true if we use the definition of $L_n$ proposed in the work [28, p.7], where the requirement $N_2(t) \in C_{-1,0}(0, \infty)$ is used.*

**Theorem 2** *Let $(M_m(t), N_m(t))$ with $m \in \mathbb{N}$ be Luchko pairs of kernels from the Luchko set $L_m$. Then the pair of the kernels*
$$M_{m,n-m}(t) = (\{1\}^{n-m} * M_m)(t), \quad N_{m,0}(t) = N_m(t), \tag{19}$$
*where $1 \leq m \leq n-1$, belongs to the Luchko set $L_n$.*

**Proof.** The condition $(M_m(t), N_m)) \in L_m$ means that
$$(M_m * N_m)(t) = \{1\}^m, \quad \text{and} \quad M_m(t), N_m(t) \in C_{-1}(0, \infty). \tag{20}$$
Then
$$(M_{m,n-m} * N_{m,0})(t) = ((\{1\}^{n-m} * M_m) * N_m)(t) =$$
$$(\{1\}^{n-m} * (M_m * N_m))(t) =$$
$$(\{1\}^{n-m} * \{1\}^m)(t) = \{1\}^n. \tag{21}$$

Therefore the kernels (19) satisfy the Luchko condition (10).

Then, using that $M_m(t), N_m(t) \in C_{-1}(0, \infty)$, we get that the kernels (19) belong to $C_{-1}(0, \infty)$. This statement is a consequence of the fact that the triple $\mathcal{R}_{-1} = (C_{-1}(0, \infty), *, +)$ is a commutative ring with the multiplication $*$ in the form of the Laplace convolution [30, 27]. Therefore the pair $(M_{m,n-m}, N_{n,0}(t))$ belongs to the Luchko set $L_n$.

□

**Definition 3** *The set of kernels pairs $(M(t), N(t))$, which belong to $L_n$ and can be represented as (19), where $(M_m(t), N_m(t))$ are the Luchko pairs of kernels from $T_m \subset L_m$ for $m \in \{1, \ldots, n-1\}$, will be denoted as $T_{n,m}$. The set $T_{n,m}$ is a subset of the Luchko set $L_n$.*

Note that as a special case of $T_{n,m}$, we can consider the kernels from the subset $T_m \subset L_m$. In general, we assume that $(M(t), N(t)) \in L_n$, which may not belong to the subset $T_m$. In other words, the kernels $M_m(t), N_m(t)$, which are used in (19) and can be represented through the Sonine pairs of kernels as (14), are only special case of $T_{n,m}$.

**Example 2** The set of kernel pairs $(M_{1,n-1}(t), N_{m,0}(t)) \in T_{n,1}$. Kernel pair (19) with $m = 1$ has the form
$$M_{1,n-1}(t) = (\{1\}^{n-1} * M_1)(t), \quad N_{1,0}(t) = N_1(t). \tag{22}$$
Note that in our case, the condition $N_1(t) \in C_{-1,0}(0, \infty)$ is not used, and we can use the pair $(M_1(t), N_1(t) \in C_{-1}(0, \infty)$ that satisfies the Sonine condition $(M_1 * N)1)(t) = \{1\}$ for $t > 0$. This set of kernel pairs $(M_{1,n-1}(t), N_{m,0}(t)) \in T_{n,1}$ was defined in [28]. For the subset $T_{n,1}$ of the Luchko set $L_n$, the general fractional derivatives and integrals were defined in [28].

**Example 3** As an example of the kernel pair from $T_{n,2}$, we can consider the kernels
$$M_{2,n-2}(t) = (\{1\}^{n-2} * M_2)(t), \quad N_{2,0}(t) = N_2(t). \tag{23}$$
As an example of the case $(M_2(t), N_2(t)) \in L_2$ such that $(M_2(t), N_2(t)) \notin T_2$, we can use
$$M_2(t) = t^{\alpha/2} J_\alpha(2\sqrt{t}), \quad N_2(t) = t^{-\alpha/2} I_{-\alpha}(2\sqrt{t}), \tag{24}$$
where $0 < \alpha < 1$, and
$$J_\alpha(t) = \sum_{k=0}^{\infty} \frac{(-1)^k (t/2)^{2k+\alpha}}{k!\,\Gamma(\alpha+k+1)}, \quad I_\alpha(t) = \sum_{k=0}^{\infty} \frac{(t/2)^{2k+\alpha}}{k!\,\Gamma(\alpha+k+1)} \tag{25}$$
with $\alpha > -1$, are the Bessel and the modified Bessel functions, respectively. The pair of kernels (24) belongs to the Luchko set $L_2$, [28, p.9].

**Example 4** As an example of the kernel pair from $T_{n,2}$, we can consider the kernels
$$M_{2,n-2}(t) = (\{1\}^{n-2} * M_2)(t), \quad N_{2,0}(t) = N_2(t), \tag{26}$$
where $(M_2(t), N_2(t)) \in T_2$. For example, we can use the kernel pairs
$$M_2(t) = (\mu_1 * \mu_2)(t), \quad N_2(t) = (\nu_1 * \nu_2)(t), \tag{27}$$
or
$$M_2(t) = (\mu_1 * \nu_2)(t), \quad N_2(t) = (\nu_1 * \mu_2)(t). \tag{28}$$
The pairs $(\mu_k(t), \nu_k(t)) \in L_1$ can be of the same or different types. For the same type, we can consider
$$\mu_k(t) = h_{\alpha_k}(t), \quad \nu_k(t) = h_{1-\alpha_k}(t), \tag{29}$$
where $\alpha_k \in (0,1)$. For the different types, we can consider $(\mu_1(t), \nu_1(t))$ given by (29), and
$$\mu_2(t) = t^{\alpha/2} J_\alpha(2\sqrt{t}), \quad \nu_2(t) = t^{-1/2-\alpha/2} I_{-\alpha-1}(2\sqrt{t}), \tag{30}$$
where $-1 < \alpha < 0$, and $J_\alpha, I_\alpha$ are given by (25). Note that pair (30) belongs to the Luchko set $L_1$, [28, p.9].

**Theorem 3** *Let $(M_m(t), N_m(t))$ be Luchko pairs of kernels from the set $T_m \subset L_m$ for $m \in \{1, \ldots, n-1\}$, and the kernels $M_{m,k}(t), N_{m,k}(t) \in C_{-1}(0, \infty)$ with $k \in \mathbb{N}$ are defined as*
$$M_{m,k}(t) = (\{1\}^k * M_m)(t), \quad N_{m,k}(t) = (\{1\}^k * N_m)(t). \tag{31}$$
*Then the pair of the kernels*
$$M_{m,n-l}(t) = (\{1\}^{n-l} * M_m)(t), \quad N_{m,l-m}(t) = (\{1\}^{l-m} * N_m)(t) \tag{32}$$
*with $m + 1 \leq l \leq n - 1$ and $1 \leq m \leq n - 1$, belongs to the Luchko set $L_n$.*

**Proof.** The condition $(M_m(t), N_m(t)) \in L_m$ means that
$$(M_m * N_m)(t) = \{1\}^m, \quad \text{and} \quad M_m(t), N_m(t) \in C_{-1}(0, \infty). \tag{33}$$
Then
$$\begin{aligned}(M_{m,n-l} * N_{m,l-m})(t) &= ((\{1\}^{n-l} * M_m) * (\{1\}^{l-m} * N_m))(t) = \\ &= (\{1\}^{n-m} * (M_m * N_m))(t) = \\ &= (\{1\}^{n-m} * \{1\}^m)(t) = \{1\}^n. \end{aligned} \tag{34}$$
Therefore the kernels (32) satisfy the Luchko condition (10).

Then, we use that $M_m(t), N_m(t) \in C_{-1}(0, \infty)$, and the fact that $\mathcal{R}_{-1} = (C_{-1}(0, \infty), *, +)$ is a commutative ring with the multiplication $*$ in the form of the Laplace convolution [30, 27]. Therefore, the kernels (32) belong to $C_{-1}(0, \infty)$.

Then using that the kernels (32) belong to $C_{-1}(0, \infty)$ and satisfy the Luchko condition (10), we state that the pair $(M_{m,n-l}, N_{m,l-m}(t))$ belongs to the Luchko set $L_n$.

□

**Definition 4** *The set of kernel pairs $(M(t), N(t))$, which belong to $L_n$ and can be represented in the form*
$$M(t) = M_{m,n-l}(t) = (\{1\}^{n-l} * M_m)(t), \quad N(t) = N_{m,l-m}(t) = (\{1\}^{l-m} * N_m)(t) \quad (35)$$
*where $(M_m(t), N_m(t))$ are the Luchko pairs of kernels from $T_m \subset L_m$ for $m \in \{1, \dots, n-1\}$, will be denoted as $T_{n,m,l}$. The sets $T_{n,m,l}$ are subsets of the Luchko set $L_n$.*

**Remark 4** *The subsets $T_{n,m,l}, T_{n,m}, T_n$, and others that are built from the Sonine pairs of kernels cannot cover the entire set of the Luchko pairs $L_n$. For example, the pair of the kernels*
$$M_n(t) = t^{\alpha/2} J_\alpha(2\sqrt{t}), \quad N_n(t) = t^{n/2-\alpha/2-1} I_{n-\alpha-2}(2\sqrt{t}), \quad (36)$$
*with $n - 2 < \alpha < n - 1$, $n \in \mathbb{N}$, belongs to the Luchko set $L_n$, [28, p.9], where $J_\alpha(t)$ and $I_\alpha(t)$ are are the Bessel and the modified Bessel functions that defined in (25).*

Using kernel pairs $(M_{m_j}(t), N_{m_j}(t))$ from the different Luchko sets $L_{m_j}$ with $j = 1, \dots p$ instead of the pairs $(\mu_j(t), \nu_j(t))$ from the Luchko set $L_1$, we can define the following generalizations of Theorems 1 and 3.

**Theorem 4** *Let $(M_{m_j}(t), N_{m_j}(t))$ with $j = 1, \dots p$ be Luchko pairs of kernels from the sets $L_{m_j}$ for $m_j \in \mathbb{N}$.*

*Then the pair of the kernels*
$$M_{\{m\},n}(t) = \left(M_{m_1} * \dots * M_{m_p}\right)(t), \quad N_{\{m\},n}(t) = \left(N_{m_1} * \dots * N_{m_p}\right)(t), \quad (37)$$
*where*
$$\sum_{j=1}^{p} m_j = n \quad (38)$$
*belongs to the Luchko set $L_n$.*

**Proof.** The proof of Theorem 4 is similar to the proof of Theorem 1. □

**Theorem 5** *Let $(M_{m_j}(t), N_{m_j}(t))$ with $j = 1, \dots p$ be Luchko pairs of kernels from the Luchko set $L_{m_j}$ for $m_j \in \mathbb{N}$.*

*Then the pair of the kernels*
$$M_{\{m\},n-l}(t) = \left(\{1\}^{n-l} * M_{\{m\},\eta}\right)(t), \quad N_{\{m\},l-\eta}(t) = \left(\{1\}^{l-\eta} * N_{\{m\},\eta}\right)(t), \quad (39)$$
*where $M_{\{m\},\eta}(t), N_{\{m\},\eta}(t)$ are defined by equations (37), and*
$$\sum_{j=1}^{p} m_j = \eta \quad (40)$$
*with $\eta + 1 \leq l \leq n - 1$ and $1 \leq \eta \leq n - 1$, belongs to the Luchko set $L_n$.*

**Proof.** The proof of Theorem 5 is similar to the proof of Theorem 3. □

The set of kernel pairs $(M(t), N(t))$, which belong to $L_n$ and can be represented by expressions that are used in Theorems 4 and 5 will be denoted as $T_{n,\{m\}}$ and $T_{n,\{m\},l}$, respectively.

Let us give some examples of the kernel pairs from the set $L_1$ (see [27, 27] and references therein).

**Example 5** The pair of the kernels
$$\mu_1(t) = \frac{t^{\alpha_1-1}}{\Gamma(\alpha_1)}, \quad (41)$$
$$\nu_1(t) = \frac{t^{-\alpha_1}}{\Gamma(1-\alpha_1)}, \quad (42)$$
is well-known in fractional calculus as kernels of the Riemann-Liouville fractional derivatives and integrals [1, 4]. This are Sonine pair of kernels, if $0 < \alpha_1 < 1$.

**Example 6** The pair of the kernels

$$\mu_2(t) = \frac{t^{\alpha_2-1}}{\Gamma(\alpha_2)} e^{-\lambda t}, \tag{43}$$

$$\nu_2(t) = \frac{t^{-\alpha_2}}{\Gamma(1-\alpha_2)} e^{-\lambda t} + \frac{\lambda^{\alpha_2}}{\Gamma(1-\alpha_2)} \gamma(1-\alpha_2, \lambda t), \tag{44}$$

and vice versa, belongs to the Sonine set, if $0 < \alpha_2 < 1$ [20, p.3627], and $\lambda \geq 0$, $t > 0$, where $\gamma(\beta, t)$ is the incomplete gamma function

$$\gamma(\beta, t) = \int_0^t \tau^{\beta-1} e^{-\tau} \, d\tau. \tag{45}$$

**Example 7** The pair of the kernels

$$\mu_3(t) = t^{\alpha_3-1} \Phi(\beta, \alpha_3; -\lambda t), \tag{46}$$

$$\nu_3(t) = \frac{\sin(\pi\alpha_3)}{\pi} t^{-\alpha_3} \Phi(-\beta, 1-\alpha_3; -\lambda t), \tag{47}$$

and vice versa belongs to the Sonine set [20, p.3629], if $0 < \alpha_3 < 1$, where $\Phi(\beta, \alpha; z)$ is the Kummer function

$$\Phi(\beta, \alpha; z) = \sum_{k=0}^{\infty} \frac{(\beta)_k}{(\alpha)_k} \frac{z^k}{k!}. \tag{48}$$

Let us give examples of the Luchko pairs of kernels from the subset $T_3 \subset L_3$. For example, the Laplace convolutions of three kernels in the from

$$M_3(t) = (\mu_1 * \mu_2 * \mu_3)(t) \quad N_3(t) = (\nu_1 * \nu_2 * \nu_3)(t), \tag{49}$$

where $\mu_k(t), \nu_k(t)$ ($k = 1,2,3$) are given in Examples 5, 6, 7, belong to the Luchko set $L_3$. This pair can be denoted as $[\mu 123, \nu 123]$, the vice versa pair as $[\nu 123, \mu 123]$. Another example is

$$M_3(t) = (\mu_1 * \mu_2 * \nu_3)(t) \quad N_3(t) = (\nu_1 * \nu_2 * \mu_3)(t) \tag{50}$$

that can be denoted as $[\mu 12\nu 3, \nu 12\mu 3]$, and vice versa pair as $[\nu 12\mu 3, \mu 12\nu 3]$.

It can be seen that there are many variations for combinations of the kernels $\mu_k(t), \nu_k(t)$ ($k = 1,2,3$), even without taking into account the possibility of using the same kernels. For example, $[\mu 13\nu 2, \nu 13\mu 2]$, $[\mu 23\nu 1, \nu 23\mu 1]$, $[\mu 1\nu 23, \nu 1\mu 23]$, $[\mu 2\nu 13, \nu 2\mu 13]$, and so on. This possibility to use different combinations of the kernels from Sonine set significantly expands the possibilities for describing nonlocalities of various types in natural and social sciences.

## 3 General Fractional Integral and Derivatives

To simplify the text, we will extend the notation of the subsets $T_{n,m,l}$, $T_{n,m}$, $T_n$, where $n \in \mathbb{N}$, $1 \leq m \leq n-1$, $m+1 \leq l \leq n-1$ by defining these notations for $l = m$ and $m = n$, respectively. Let us define them by the expressions

$$T_{n,m,m} := T_{n,m}, \quad T_{n,n} := T_n. \tag{51}$$

Using these notation, we give the definitions of the general fractional operator for the pairs of kernels from $T_{n,m,m}$ with $1 \leq m \leq n$, and $m \leq l \leq n$.

**Definition 5** Let $(M_{m,n-l}(t), N_{m,l-m}(t))$ with $1 \leq m \leq n$, and $m \leq l \leq n$ be the Luchko pair of kernels from the subset $T_{n,m,l} \subset L_n$. Then, the GF-integral is defined by the equation

$$I_{(M)} X(t) = (M_{m,n-l} * X)(t) = \int_0^t M_{m,n-l}(t-\tau) X(\tau) \, d\tau \tag{52}$$

for $X(t) \in C_{-1}(0, \infty)$.

**Definition 6** Let $(M_{m,n-l}(t), N_{m,l-m}(t))$ $1 \leq m \leq n$, and $m \leq l \leq n$ be the Luchko pair of kernels from the subset $T_{n,m,l} \subset L_n$. Then, the GF-derivative of Caputo type is defined by the equation

$$D^*_{(N)} X(t) = (N_{m,l-m} * X^{(n)})(t) = \int_0^t N_{m,l-m}(t-\tau) X^{(n)}(\tau) \, d\tau \tag{53}$$

for $X(t) \in C^n_{-1}(0, \infty)$, where

$$C^n_{-1}(0, \infty) := \{X: \ X^{(n)}(t) \in C_{-1}(0, \infty)\}. \tag{54}$$

**Definition 7** Let $(M_{m,n-l}(t), N_{m,l-m}(t))$ $1 \leq m \leq n$, and $m \leq l \leq n$ be the Luchko pair of kernels from the subset $T_{n,m,l} \subset L_n$. Then, the GF-derivative of Riemann-Liouville type is defined by the equation

$$D_{(N)} X(t) = \frac{d^n}{dt^n}(N_{m,l-m} * X)(t) = \frac{d^n}{dt^n} \int_0^t N_{m,l-m}(t-\tau) X(\tau) \, d\tau \tag{55}$$

for $X(t) \in C_{-1}(0, \infty)$.

The GF operators are defined similarly for subsets $T_{n,m}$ and $T_{n,m}$ of the Luchko set $L_n$. For the subsets $T_{n,\{m\}}$ and $T_{n,\{m\},l}$ of $L_n$, the GFI and GFDs are are defined similarly to Definitions 5, 6, 7.

## 4 Fundamental Theorems of General Fractional Calculus

The fundamental theorems (FT) of standard calculus for derivatives and integrals of integer order $n \in \mathbb{N}$ are the following. The first FT is written as

$$\frac{d^n}{dt^n} I^n X(t) = X(t). \tag{56}$$

The second FT is given as

$$I^n X^{(n)}(t) = X(t) - \sum_{k=0}^{n-1} X^{(k)}(0) h_{k+1}(t). \tag{57}$$

Here $I^n$ is the integral of the order $n \in \mathbb{N}$ such that

$$I^n X(t) = \int_0^t d\tau_1 \int_0^{\tau_1} d\tau_2 \dots \int_0^{\tau_{n-1}} d\tau_n X(\tau_n) =$$
$$\frac{1}{(n-1)!} \int_0^t (t-\tau)^{n-1} X(\tau)\, d\tau = \int_0^t h_n(t-\tau) X(\tau)\, d\tau. \tag{58}$$

Let us prove FT of the general fractional calculus for multi-kernel approach.

**Theorem 6 (First FT for the GF-derivative of Caputo type)**

*Let $(M_{m,n-l}(t), N_{m,l-m}(t))$ with $1 \leq m \leq n-1$, and $m+1 \leq l \leq n-1$ be the Luchko pair of kernels from the subset $T_{n,m,l} \subset L_n$.*

*Then, for the GF-derivative $D^*_{(M)}$ of Caputo type (6) and the GF-integral (5), the equalities*

$$D^*_{(N)} I_{(M)} X(t) = X(t) \tag{59}$$

*holds for $X(t) \in C_{-1,(N)}(0, \infty)$, where*

$$C_{-1,(N)}(0, \infty) := \{X: X(t) = I_{(N)} Y(t) = (N_{m,l-m} * Y)(t),\ Y(t) \in C_{-1}(0, \infty)\}. \tag{60}$$

**Proof.** The condition $X(t) \in C_{-1,(N)}(0, \infty)$ means that

$$X(t) = I_{(N)} Y(t) = (N_{m,l-m} * Y)(t). \tag{61}$$

The property $(M_{m,n-l}(t), N_{m,l-m}(t)) \in L_n$ leads to

$$(M_{m,n-l} * N_{m,l-m})(t) = \{1\}^n. \tag{62}$$

Using these equalities, we get

$$D^*_{(N)} I_{(M)} X(t) = D^*_{(N)} I_{(M)} I_{(N)} Y(t) =$$
$$D^*_{(N)} \left( M_{m,n-l} * (N_{m,l-m} * Y) \right)(t) =$$
$$D^*_{(N)} \left( (M_{m,n-l} * N_{m,l-m}) * Y \right)(t) = D^*_{(N)} (\{1\}^n * Y)(t).$$

Using the definition of the operators $D^*_{(N)}$ and property (61), we obtain

$$D^*_{(N)} (\{1\}^n * Y)(t) = D^*_{(N)} I^n Y(t) =$$
$$\left( N_{m,l-m} * \frac{d^n}{dt^n}(I^n Y) \right)(t) =$$
$$(N_{m,l-m} * Y)(t) = I_{(N)} Y(t) = X(t), \tag{63}$$

where $I^n$ is defined by (58). □

**Theorem 7 (Second FT for the GF-derivative of Caputo type)** *Let $(M_{m,n-l}(t), N_{m,l-m}(t))$ with $1 \leq m \leq n-1$, and $m+1 \leq l \leq n-1$ be the Luchko pair of kernels from the subset $T_{n,m,l} \subset L_n$.*

*Then, for the GF-derivative $D^*_{(M)}$ of Caputo type (6) and the GF-integral (5), the equality*

$$I_{(M)} D^*_{(N)} X(t) = X(t) - \sum_{k=0}^{n-1} X^{(k)}(0) h_{k+1}(t) \tag{64}$$

*holds for $X(t) \in C^n_{-1}(0, \infty)$, where*

$$C^n_{-1}(0, \infty) := \{X: X^{(n)}(t) \in C_{-1}(0, \infty)\}. \tag{65}$$

**Proof.** The property $(M_{m,n-l}(t), N_{m,l-m}(t)) \in L_n$ means that

$$(M_{m,n-l} * N_{m,l-m})(t) = \{1\}^n. \tag{66}$$

Using the definitions of operators $I_{(M)}$, $D^*_{(N)}$, and the fundamental theorem of standard calculus

$$I^n X^{(n)}(t) = X(t) - \sum_{k=0}^{n-1} X^{(k)}(0) h_{k+1}(t), \tag{67}$$

we obtain

$$I_{(M)} D^*_{(N)} X(t) = \left(M_{m,n-l} * \left(N_{m,l-m} * X^{(n)}\right)\right)(t) =$$
$$\left(\left(M_{m,n-l} * N_{m,l-m}\right) * X^{(n)}\right)(t) =$$
$$(\{1\}^n * X^{(n)})(t) = I^n X^{(n)}(t) =$$
$$X(t) - \sum_{k=0}^{n-1} X^{(k)}(0) h_{k+1}(t). \tag{68}$$

□

**Theorem 8 (First FT for the GF-derivative of RL type)**

Let $(M_{m,n-l}(t), N_{m,l-m}(t))$ with $1 \leq m \leq n-1$, and $m+1 \leq l \leq n-1$ be the Luchko pair of kernels from the subset $T_{n,m,l} \subset L_n$.

Then, for the GF-derivative $D_{(M)}$ of Riemann-Liouville type (7) and the GF-integral (5), the equalities
$$D_{(N)} I_{(M)} X(t) = X(t) \tag{69}$$
holds for $X(t) \in C_{-1}(0, \infty)$.

**Proof.** Using the definitions of the GF-derivative $D_{(M)}$ of Riemann-Liouville type (7) and the GF-integral (5), and the standard fundamental theorem
$$\frac{d^n}{dt^n} I^n X(t) = X(t), \tag{70}$$
we obtain
$$D_{(N)} I_{(M)} X(t) = \frac{d^n}{dt^n} \left(N_{m,l-m} * \left(M_{m,n-l} * X\right)\right)(t) =$$
$$\frac{d^n}{dt^n} \left(\left(N_{m,l-m} * M_{m,n-l}\right) * X\right)(t) = \frac{d^n}{dt^n}(\{1\}^n * X)(t) =$$
$$\frac{d^n}{dt^n} I^n X(t) = X(t). \tag{71}$$

□

**Theorem 9 (Second FT for the GF-derivative of RL type)** Let $(M_{m,n-l}(t), N_{m,l-m}(t))$ with $1 \leq m \leq n-1$, and $m+1 \leq l \leq n-1$ be the Luchko pair of kernels from the subset $T_{n,m,l} \subset L_n$.

Then, for the GF-derivative $D_{(M)}$ of Riemann-Liouville type (7) and the GF-integral (5), the equalities
$$I_{(M)} D_{(N)} X(t) = X(t). \tag{72}$$
holds for $X(t) \in C_{-1,(M)}(0, \infty)$, where
$$C_{-1,(M)}(0, \infty) := \{X: \ X(t) = I_{(M)} Y(t) = (M_{m,n-l} * Y)(t), \ Y(t) \in C_{-1}(0, \infty)\}. \tag{73}$$

**Proof.** The condition $X(t) \in C_{-1,(M)}(0, \infty)$ means that
$$X(t) = I_{(M)} Y(t) = (M_{m,n-l} * Y)(t). \tag{74}$$
The property $(M_{m,n-l}(t), N_{m,l-m}(t)) \in L_n$ leads to
$$(M_{m,n-l} * N_{m,l-m})(t) = \{1\}^n. \tag{75}$$
Using these equalities and definition of the operator $D_{(N)}$, we get
$$I_{(M)} D_{(N)} X(t) = I_{(M)} D_{(N)} I_{(M)} Y(t) =$$
$$I_{(M)} \frac{d^n}{dt^n} \left(N_{m,l-m} * M_{m,n-l}(t) * Y\right)(t) =$$
$$I_{(M)} \frac{d^n}{dt^n} (\{1\}^n * Y)(t) = I_{(M)} \frac{d^n}{dt^n} I^n Y(t) = I_{(M)} Y(t) = X(t). \tag{76}$$

□

The Fundamental Theorems of GF derivatives and GF integrals are proves similarly for subsets $T_{n,m}$ and $T_n$ of the Luchko set $L_n$.

For the subsets $T_{n,\{m\}}$ and $T_{n,\{m\},l}$ of $L_n$, the Fundamental Theorems for the GFI and GFDs are proves similarly to proofs of Theorems 6, 7, 8, 9.

## 5 Conclusion

In this work, we develop approaches to building theory of general fractional calculus of arbitrary order, which are proposed in the works [27, 28]. We propose an extension of Luchko approach to formulate the general fractional calculus of arbitrary order for multi-kernel cases. In the suggested approach, we use different types (subsets) of the pairs of kernels in definitions of the general fractional integrals and derivatives. The proposed set of kernel pairs are new subsets of the Luchko set of kernel pairs, which were not considered in the works [27, 28]. The first and second fundamental theorems for the proposed general fractional derivatives and integrals are proved.

In a more general form, the main idea of this article is as follows: We propose to use different generators (generating set) of pairs of kernels for each Luchko set $L_k$ and powers $\{1\}^k$ with $k = 1, \dots, n-1$ to construct pairs of kernels for Luchko set $L_n$ of the next order $n$. The possibility of this approached is based on the fact that the triple $\mathcal{R}_{-1} = (C_{-1}(0,\infty), *, +)$, where the multiplication $*$ is the Laplace convolution and $+$ the standard addition of functions, is a commutative ring without divisors of zero [30, 27]. Here we mean that a generating set of the Luchko set $L_k$ is a subset of $L_k$ such that almost every kernel of $L_k$ can be expressed as a combination (by using the Laplace convolution) of finitely many kernels of the subset and their associated kernels (their inverses). If the ring $\mathcal{R}_{-1}$ is a ring having a system of generators $G_n$, then almost every kernel can be represented as a product (the Laplace convolution) of kernels from $G_n$, and inverse (associated) to them. Note that the number of kernels multiplied by the Laplace convolution can be more than $n$ (for example, $(h_\alpha * h_\beta)(t) = h_{\alpha+\beta}(t)$). The mathematical implementation of this idea is a complex and interesting problem that should to be solved. It is not obvious that this idea can be implemented, but research in this direction can lead to interesting results, which will be important for general fractional calculus and its various applications.

The proposed approach to general fractional calculus can be useful for various applications in physics, economics, nonlinear dynamics and for other areas of science. In applications of the general fractional calculus, it is useful to have a wider range of operators, for example, which kernels are defined by using the Laplace convolution of different types of kernels. The importance of the proposed approach to GFC is related with the importance of describing systems and processes with a wider variety of non-localities in time and space [8, 9, 10, 11, 12], [31, 32, 33]. The GFC and the proposed multi-kernel approach to GFC can be important to obtaining results concerning of general form of nonlocality, which can be described by general form operator kernels, and not its particular implementations and representations [34]. For example, we can derive general nonlocal maps as exact solutions of nonlinear equations with GFI and GFD at discrete points [34], and a general approach to the construction of non-Markovian quantum theory can be proposed.